\documentclass[12pt]{article}
\usepackage[english]{babel}
\usepackage{graphicx} 
\usepackage{latexsym,color,amsmath,amsthm,amssymb,amscd,amsfonts}
\usepackage{scalefnt}
\usepackage[font=small,labelfont=bf]{caption}
\usepackage{float}
\usepackage{graphicx}
\usepackage{amsmath}
\usepackage{relsize}
\usepackage{mathrsfs}
\usepackage{upgreek}
\usepackage{sectsty}
\usepackage{bbm}

\newtheoremstyle{nonum}{}{}{\itshape}{}{\bfseries}{.}{ }{\thmnote{#3}}

\makeatletter
\newcommand{\cbigoplus}{\DOTSB\cbigoplus@\slimits@}
\newcommand{\cbigoplus@}{\mathop{\widehat{\bigoplus}}}
\makeatother

\newtheorem{thm}{Theorem}[section]
\newtheorem*{thm*}{Theorem}

\newtheorem{cor}[thm]{Corollary}
\newtheorem{lem}[thm]{Lemma}
\newtheorem*{lem*}{Lemma}

\newtheorem*{rem*}{Remark}
\newtheorem{prop}[thm]{Proposition}
\newtheorem*{prop*}{Proposition}

\newtheorem*{definition*}{Definition}

\newtheorem*{fact*}{Fact}

\newtheorem*{rems*}{Remarks}

\newtheorem{conj}[thm]{Conjecture}
\theoremstyle{nonum}

\newcommand{\R}{\mathbb R}

\def\Vol{{\rm Vol}}

\usepackage{bbm}

\def\eps{{\varepsilon}}

\def\conv{{\rm conv}}

\def\Vol{{\rm Vol}}
\def\diam{{\rm diam}}

\title{On the Figiel-Lindenstrauss-Milman inequality}
\author{Tomer Milo}
\date{}

\begin{document}

\maketitle
\begin{abstract}
\noindent The Figiel-Lindenstrauss-Milman inequality is a fundamental inequality in the combinatorial theory of polytopes. It is classically obtained as a corollary of Milman's version of Dvoretzky's theorem. The goal of this paper is to provide a short and elementary proof of this inequality, derive more general versions of it, and discuss its tightness, where much is not known.
\end{abstract}

\section{Introduction}

The Figiel-Lindenstrauss-Milman (FLM) inequality is the following inequality, originally obtained in \cite{FLM} as an application Milman's version of Dvoretzky's Theorem.

\begin{thm}\label{Theorem: FLM}(FLM inequality) There exists an absolute constant $c>0$ such that for every $n \in \mathbb{N}$ and for every convex polytope $P \subset \mathbb{R}^{n}$ satisfying $rB_2^n \subset P \subset RB_2^n$:
\[ \log |V| \cdot \log|\mathcal{F}| \cdot \left( \frac{R}{r} \right)^2 \geq cn^2 \]
Where $V$ and $\mathcal{F}$ are the sets of vertices ($0$-dimensional) and facets ($(n-1)$-dimensional facets) of $P$, respectively, and $B_2^n$ denotes the Euclidean unit ball in $\mathbb{R}^n$ centered at the origin.
\end{thm}

\noindent This theorem has the following immediate corollary, which is sometimes also referred to as the FLM inequality in the literature. We shall refer to it throughout as the `symmetric FLM inequality'.

\begin{cor}\label{cor: FLM corollary}(Symmetric FLM inequality)
    There exists an absolute constant $c>0$ such that for every $n \in \mathbb{N}$ and for every convex polytope $P \subset \mathbb{R}^{n}$ which is centrally symmetric (meaning $P = -P$):
    \[ \log|V| \cdot \log|\mathcal{F}| \geq cn \]
\end{cor}

\noindent Corollary \ref{cor: FLM corollary} follows from Theorem \ref{Theorem: FLM} by using John's Theorem (see \cite[chapter 2]{AGA} for a detailed exposition), which states that if $P$ is centrally symmetric, then one can find an invertible linear transformation $T$ such that $B_2^n \subset T(P) \subset \sqrt{n}B_2^n$. Since invertible linear transformations do not change the combinatorial structure of a polytope, the corollary is proven. The goal of this note is to provide a short and elementary proof of Theorem \ref{Theorem: FLM} and discuss its tightness. As we shall see, the tightness of the symmetric FLM inequality (or the case $\frac{R}{r} = \sqrt{n}$) is completely understood, while for any other choice of growth for $\frac{R}{r}$ the tightness question is essentially completely open. See section 3 for further details. 

\subsection{Preliminaries and notations}

We call a subset $P \subset \mathbb{R}^n$ a \textbf{polytope} if it is the convex hull of finitely many points and has a non-empty interior. We denote by $V_P$ and $\mathcal{F}_P$ the sets of vertices ($0$-dimensional faces) and facets ($(n-1)$-dimensional faces) respectively, and often just write $V$ and $\mathcal{F}$ when the polytope which we refer to is clear from the context. When we have a sequence of polytopes $P_n \subset \mathbb{R}^n$ we will also often write $V_n$ instead of $V_{P_n}$, etc. For two quantities $A$ and $B$, which may depend on the dimension, a polytope, and other parameters, we say that $A \lesssim B$ if there exists some absolute constant $c>0$ so that $A < cB$. We write $A \sim B$ if $A \lesssim B$ and  $B \lesssim A$. The letter $c$ denotes some positive universal constant, unless specifically stated otherwise.

\noindent For a polytope $P \subset \mathbb{R}^n$, we denote by $P^{\circ} := \{x\in \mathbb{R}^n : \max_{y\in P} \langle x,y \rangle \leq 1\}$ the \textbf{dual polytope} of $P$ (this is indeed a polytope whenever $0 \in {\rm int}(P)$). It is easy to verify that the dual operation switches between the numbers of vertices and facets of a polytope; that is, $|V_{P^{\circ}}| = |\mathcal{F}_P|$ and $|\mathcal{F}_{P^{\circ}}| = |V_P|$, and that duality is an involution: $P^{\circ \circ} = P$. The \textbf{norm} of $P$ is given by $\|x\|_P = \inf\{t: x \in tP\}$ (this is only a true norm for a centrally symmetric $P$, but we shall abuse the language and always refer to it as a norm). Moreover, the dual norm is given by: $\| x \|_{P^{\circ}} = \max_{y\in P} \langle x, y \rangle$, and $\|x\|_{P^{\circ}}$ is referred to as the \textbf{support function} of $P$ and is denoted by $h_{P}$(x). As is standard, we shall use the following notation:
\[ M(P) := \int_{S^{n-1}} \|\theta\|_{P}d\sigma(\theta),  \quad  M^{*}(P) := \int_{S^{n-1}}h_{P}(\theta)d\sigma \]
Where $\sigma$ denotes the unique Haar probability measure on $S^{n-1}$. Note that $M^{*}(P) = M(P^{\circ})$.

\noindent We will often use the following elementary facts: given an orthogonal decomposition $\mathbb{R}^{n+k} = \mathbb{R}^n \bigoplus \mathbb{R}^{k}$ and polytopes $P \subset \mathbb{R}^n$ and $Q \subset \mathbb{R}^k$, then for the polytope $P \times Q \subset \mathbb{R}^{n+k} = \mathbb{R}^n \bigoplus \mathbb{R}^k$, we have  $|V_{P \times Q}| = |V_P||V_Q|$ and $|\mathcal{F}_{P \times Q}| = |\mathcal{F}_P| + |\mathcal{F}_Q|$. Finally, it is elementary to verify that $(P \times Q)^{\circ} = \conv(P^{\circ}, Q^{\circ})$, where the duals $P^\circ$ and $Q^\circ$ are computed in the respective subspaces $\mathbb{R}^n$ and $\mathbb{R}^k$ of $\mathbb{R}^{n+k}$. All logarithms are to base $e$ unless specifically stated otherwise.

\section{Proof of the FLM inequality}

The goal of this section is to prove Theorem \ref{Theorem: FLM}. It is classically proved as a consequence of Dvoretzky's theorem; see \cite[section 5]{AGA} for an exposition of Milman's classical proof and its applications. Here we provide a new proof avoiding Dvoretzky's theorem completely (as well as the Dvoretzky-Rogers lemma), and also avoiding the classical concentration inequality for Lipschitz functions on the Euclidean sphere $S^{n-1} \subset \mathbb{R}^n$, which is proved via the spherical isoperimetric inequality due to Levy. We shall only need the following simple concentration result, for which an elementary proof can be found in \cite[Theorem 3.1.5]{AGA}:

\begin{lem}\label{lemma: concentration of measure} Let $\sigma$ denote the unique rotational-invariant probability measure on the Euclidean unit sphere $S^{n-1}$. For every $u \in S^{n-1}$ and every $\eps \in (0,1)$,

\[ \sigma(\{\theta \in S^{n-1}: \langle \theta, u \rangle < \eps \}) > 1 - e^{-\frac{n\eps^2}{2}}. \]

\end{lem}
\noindent We prove Theorem \ref{Theorem: FLM} for a general polytope $ P = \conv(V) \subset \mathbb{R}^n$ satisfying $rB_n \subset P \subset R B_n$. The following lemma is the main ingredient in the proof:

\begin{lem}\label{lemma: mean width bound}
    Let $P \subset RB_n$ be a polytope with vertex set $V$ such that $\log |V| < \frac{n}{3}$. Then
    \[ M^*(P) := \int_{S^{n-1}} h_{P}d\sigma \leq CR \sqrt{\frac{\log |V|}{n}} \]
    Where $C = C_n = \sqrt{3} + o(1)$ as $n \to \infty$.
\end{lem}
\noindent \textbf{Remark.} This inequality is essentially the classical estimate on the expected value of the maximum of $|V|$ Gaussians with mean $0$ and standard deviation $R$.
\begin{proof}
We write $h_{P}(\theta) = \max_{v \in V} \langle v, \theta \rangle$. Let $B_t = \{\theta \in S^{n-1}: h_{P}(\theta) \leq t \}$. Using Lemma \ref{lemma: concentration of measure of measure}  and a union bound again gives:
    \[ M^*(P) = \int_{S^{n-1}}\max_{v \in V} \langle v, \theta \rangle d\sigma = \int_{B_t}\max_{v \in V} \langle v, \theta \rangle d\sigma  + \int_{S^{n-1} \setminus B_t }\max_{v \in V} \langle v, \theta \rangle d\sigma \leq t + R|V|e^{-\frac{1}{2}n(\frac{t}{R})^2}.
    \] 
Choosing $t = R\sqrt{\frac{3\log|V|}{n}}$ and using the fact that $|V| > n $ gives
\[ M^{*}(P) \leq  R\left( \sqrt{\frac{3\log|V|}{n}} + \frac{1}{\sqrt{|V|}}\right) \leq \left(\sqrt{3} + o(1)\right) R \sqrt{\frac{\log|V|}{n}}, \]

\noindent where the $o(1)$ term has a decay rate of at most $\frac{1}{\sqrt{\log n}}$.

\end{proof}

\noindent We apply Lemma \ref{lemma: mean width bound} to the dual polytope $P^{\circ}$. We use the following facts: $M^*(P) = M(P^\circ)$, $R(P) = \frac{1}{r(P^{\circ})}$ where $r = r(P^\circ)$ is the maximal $r$ for which $rB_n \subset P^{\circ}$, and that the number of facets $\mathcal{F}$ of $P$ is the same as the number of vertices of $P^{\circ}$. We get:
\begin{equation}
 M(P) \leq \left(\sqrt{3} + o(1)\right)\frac{1}{r}\sqrt{\frac{\log|\mathcal{F}|}{n}}.    
\end{equation}

\begin{proof}[Proof of Theorem \ref{Theorem: FLM}]
Combining Lemma \ref{lemma: mean width bound} with (1) (note that if the condition $\log|V| < \frac{n}{3}$ of Lemma \ref{lemma: mean width bound} is not satisfied, there is nothing to prove) with the trivial bound $M(P)M^*(P) \geq 1$ we get:

\[ \log|V| \cdot \log|\mathcal{F}| \geq \left(\frac{1}{9} + o(1)\right) n^2 \left( \frac{r}{R} \right)^2 (M(P)M^*(P))^2 \geq \left(\frac{1}{9} + o(1)\right) \left( \frac{r}{R} \right)^2 n^2. \]

\end{proof}

\section{The refined FLM inequality}
The following was proved in \cite{Bar} by Barvinok. 
\begin{thm}\label{Theorem: barvinok} For any $\alpha > 0$ there exists $\gamma >0$ such that the following holds: for any dimension $n$ and any centrally symmetric polytope $P \subset \mathbb{R}^n$, if $\mathcal|{\mathcal{F}|} \leq \alpha n$ then $|V| \geq e^{\gamma n}$. The same result holds when one assumes that $B_2^n \subset P \subset \sqrt{n}B_2^n$ instead of the central-symmetry assumption.
\end{thm}
\noindent Of course, by duality, the same holds when the roles are reversed: linearly many vertices implies exponentially many vertices, whenever $\frac{R}{r} \leq \sqrt{n}$. In this case, the symmetric FLM inequality becomes: $\log|V| \cdot \log|\mathcal{F}| \gtrsim n\log n$, where the implied constant depends only on $\alpha$.

\noindent The first goal of this section is to provide a generalization of Theorem \ref{Theorem: barvinok}. We will prove that a bound of the form $|\mathcal{F}|\leq \alpha n$ implies that $P$ has very high-dimensional projections which are close to Euclidean (which automatically implies that $P$ has many vertices). This will be obtained as a corollary of the following result of Gluskin (see \cite{Glus}).
\begin{thm}\label{Gluskin} Let $x_1,...,x_N \in S^{n-1}$. Let $P_N$ denote their convex hull. Then:
\[ \left( \frac{\Vol(P_N)}{\Vol(B_2^n)}\right)^{\frac{1}{N}} \lesssim \sqrt{\frac{\log(1+\frac{N}{n})}{n}} \]
\end{thm}
\begin{cor}\label{cor: refined FLM}(The refined FLM inequality)
    Let $P \subset \mathbb{R}^n$ be a polytope such that $B_2^n \subset P \subset \sqrt{n}B_2^n$. Suppose $|\mathcal{F}| = \alpha n$. Then:
    \[ \log|V| \cdot \log|\mathcal{F}| \gtrsim n\cdot \frac{\log (\alpha n)}{\log (1+\alpha)} \]
\end{cor}
\noindent \textbf{Remark.} This inequality is asymptotically stronger than the symmetric FLM inequality whenever $\alpha = \alpha_n = o(n^{\eps})$ for all (fixed) $\eps >0$.

\begin{proof}
By Holder's inequality and Theorem \ref{Gluskin} (note that $P^{\circ} \subset B_2^n$):
    \[ M^{*}(P) = M(P^{\circ}) = \int_{S^{n-1}}\|x\|_P^{\circ}d\sigma \geq \left( \int_{S^{n-1}}\|x\|_{P^{\circ}}^{-n}d\sigma \right)^{-\frac{1}{n}} = \left( \frac{\Vol(B_2^n)}{\Vol(P^{\circ})} \right)^{\frac{1}{n}} \gtrsim \sqrt{\frac{n}{\log \left(1+\alpha \right)}} \]
By a variant of Milman's version of Dvoretzky's theorem (see \cite[chapter 5]{AGA} for an exposition), there exists a subspace $E \subset \mathbb{R}^n$ of dimension $k \sim \frac{n}{\log \left(1+\alpha \right)} $ such that (after applying some invertible linear transformation):
\[ E \cap B_2^n \subset Proj_{E}(P) \subset 2(E \cap B_2^n) \]
By Theorem \ref{thm: keith ball bound} below, this means that the polytope $Proj_{E}(P)$ has $\gtrsim e^{\frac{cn}{\log \left(1+\alpha \right)}} $ vertices. Since projections can only decrease the number of vertices of the original polytope, we have
\[ \log|V| \gtrsim \frac{n}{\log (1+\alpha)}, \]
and finally:
\[ \log|V| \cdot \log|\mathcal{F}| \gtrsim n \cdot \frac{\log(\alpha n)}{\log (1+\alpha)}. \]
\end{proof}
\noindent \textbf{Remarks.}1. The same argument shows that for $P \subset \mathbb{R}^n$ such that $rB_2^n \subset P \subset RB_2^n$,
\begin{equation} \log|V| \cdot \log|\mathcal{F}| \cdot \left( \frac{R}{r} \right)^2 \gtrsim n^2\cdot \frac{\log (\alpha n)}{\log(1+\alpha)} 
\end{equation}
2. This type of argument (polytopes having many vertices/faces because of some high dimensional section/projection which is almost Euclidean) will be used many times throughout this paper.

\section{Tightness of the FLM inequality}
\subsection{The symmetric FLM inequality}

We begin the discussion of the tightness of the FLM inequality with the symmetric case, where essentially everything is known.

\begin{prop}\label{prop: tightness with a}
    For every number $a \in (0,1)$ There exists polytopes $P_n \subset \mathbb{R}^n$ such that 
    \[ \log|V_n| \sim n^a \quad and \quad \log|\mathcal{F}_n| \sim n^{1-a} \]
\end{prop}

\noindent \textbf{Remark.} The case $a= \frac{1}{2}$ was treated in \cite{FLM}. We reproduce their construction here and generalize it.

\begin{proof}
We start by considering the case $a=\frac{1}{2}$. We construct polytopes $P_{2^n} \subset \mathbb{R}^{2^n}$ such that $\log|V_{2^n}| \sim \log|\mathcal{F}_{2^n}| \sim 2^{\frac{n}{2}}$, by taking Cartesian products and convex hulls iteratively to obtain a certain Hanner polytopes, as follows: Starting with \\ $P_{0} = [-1,1] \subset \mathbb{R}$, and assuming $P_{2^n} \subset \mathbb{R}^{2^n}$ is defined, we define $P_{2^{n+1}} \subset \mathbb{R}^{2^{n+1}} = \mathbb{R}^{2^n} \bigoplus \mathbb{R}^{2^n}$ as
\[ P_{2^{n+1}} = \begin{cases} \conv(P_{2^n}, P_{2^n}), \: n \: even \\
P_{2^{n}} \times P_{2^{n}}, \: n \: odd \\
\end{cases}
\]
Clearly $P_{2^n}$ is centrally symmetric for all $n$, and 
\[ |V_{2^{n+1}}| = \begin{cases} 
2|V_{2^n}|, \: n \: even \\ 
|V_{2^n}|^2,\:  n \: odd \\
\end{cases} \]
and by duality,
\[ |\mathcal{F}_{2^{n+1}}| = \begin{cases} 
2|\mathcal{F}_{2^n}|, \: n \: odd \\ 
|\mathcal{F}_{2^n}|^2,\:  n \: even \\ 
\end{cases} \]
It is now obvious that we have the following asymptotics for $P_n$: 
\[ \log|V_{2^n}| \sim \log|\mathcal{F}_{2^n}| \sim 2^{\frac{n}{2}}. \]
which is what we wanted to prove. To extend the construction to all dimensions, we proceed as follows: fix a dimension $n$ and write $n = \sum_{i\in I}2^{r_i}$ for some index set $I \subset [k]$ and some natural numbers $r_1 >..>r_{|I|}$, where $r_1 = k = \lfloor \log n\rfloor$. Orthogonally decomposing $\mathbb{R}^n = \bigoplus_{i \in I} \mathbb{R}^{2^{r_i}}$ and taking $P_{2^{r_i}} \subset \mathbb{R}^{2^{r_i}}$ be the polytopes previously constructed, we let
\[ P_n := P_{2^{r_1}} \times ... \times P_{2^{r_{|I|}}} \subset \mathbb{R}^n,\]
we have 
\[ \log|V_n| = \sum_{i\in I} 2^{\frac{r_i}{2}} \sim 2^{\frac{r_1}{2}} \sim 2^{\frac{\log n}{2}} =\sqrt{n} \] 
and similarly,  
\[ \log|\mathcal{F}_n| = \log \left( \sum_{i \in I}|\mathcal{F}_{2^{r_i}}| \right) \sim 2^{\frac{r_1}{2}} \sim 2^{\frac{\log n}{2}} = \sqrt{n} . \] 

\noindent To generalize this argument to any $a \in (0,1)$, we take products `$a$ of the time' instead of half of the time; more precisely, let $A \subset \mathbb{N}$ be any set of density $a$ (say, $A = \{ \lfloor \frac{n}{a} \rfloor: n \in \mathbb
{N} \}$). Let
    \[ P_{2^{n+1}} = \begin{cases} \conv(P_{2^n}, P_{2^n}), \quad n \in A \\
P_{2^{n}} \times P_{2^{n}}, \quad else \\
\end{cases}\]

\noindent So $P_{2^n} \subset \mathbb{R}^{2^n}$, and an analogous argument to the previous one shows that $\log|V_n| \sim 2^{an}$ and $\log|\mathcal{F}_n| \sim 2^{(1-a)n}$. As before, by decomposing a general $n \in \mathbb{N}$ as $\sum_{i\in I}2^{r_i}$ for some index set $I \subset [k]$ and some natural numbers $r_1 >..>r_{|I|}$, where $r_1 = k \leq \lfloor \log_{\lfloor \frac{1}{a} \rfloor} n \rfloor$, and taking products of the polytopes we have already constructed ($P_n = \bigoplus_{i \in I} P_{2^{r_i}} $) we expand our construction to all dimensions. Note that the implied constants in the asymptotics of $\log|V_n|$ and $\log|\mathcal{F}_n|$ depend only on $a$.
\end{proof}
\noindent Now that we have dealt with the case $\log|V_n| \sim n^a$ and $\log|\mathcal{F}_n| \sim n^{1-a}$, one can ask more generally: Given a `growth function' $f(n) \to_{n \to \infty} \infty$, when can we find polytopes $P_n \subset \mathbb{R}^n$ such that $\log|V_n| \sim \frac{n}{f(n)}$ and $\log |\mathcal{F}_n| \sim f(n)$? The following proposition makes this precise.

\begin{prop}\label{prop: general tightness}
Let $f(n) \to_{n \to \infty} \infty$ be some arbitrary sequence of integers such that $f(n) < \sqrt{\frac{n}{2}}$. There exists polytopes $P_N \subset \mathbb{R}^{N}$, where $N = N(n)$ satisfies $\frac{n}{2} \leq N \leq n$, such that:
\[ |V_N| \sim 2^{c\frac{n}{f(n)}} \: \: and \: \: |\mathcal{F}_N| \sim n \frac{2^{c'f(n)}}{f(n)}, \]
where $c,c'>0$ are absolute constants. In particular,
\begin{equation}
    \log|V_N| \cdot \log|\mathcal{F}_N| \sim N\left( 1 + \frac{\log N}{f(N)} \right) . 
\end{equation} 

\end{prop}

\noindent \textbf{Remark.}
Choosing $f(N) = \log N $ gives polytopes $P_N \subset \mathbb{R}^N$ with $\log|\mathcal{F}_N| \sim \log N$ and $\log|V_N| \sim \frac{N}{\log N}$, showing that Barvinok's Theorem \ref{Theorem: barvinok} cannot be extended to the case where the number of facets has a polynomial upper bound rather than a linear one.
\\

\noindent Proposition \ref{prop: general tightness} demonstrates that the symmetric FLM inequality is tight, up to a constant, for almost all possible choices of sizes of $|V_N|$ and $|\mathcal{F}_N|$, in the following sense: whenever $f(N)$ has at least a logarithmic type growth (that is, there exists some $\eps >0$ such that $f(N) > \eps \log N$), we have $\log|V_N| \cdot \log|\mathcal{F}_N| \sim N$ in (3) (where the implied constant depends only on $\eps$). Noting that $f(N) = \eps \log N$ corresponds to $|\mathcal{F}_N| \sim \frac{N^{1+c'\eps}}{\eps^2\log^2N}$, we can conclude that whenever a sequence $a_n$ has a type of growth of at least $n^{1+\eps}$ for some $\eps >0$ (and also $\log(a_n) < \sqrt{n}$) one can generate polytopes $P_N \subset \mathbb{R}^N$ for which $\log|\mathcal{F}_N| \sim \log a_N$ and $\log|V_N| \sim \frac{N}{\log a_N}$, by writing $a_N = N \frac{2^{f(N)}}{cf(N)^2}$ and inverting the relation to find the appropriate $f$ to apply Proposition \ref{prop: general tightness} to (note that $f(n) \sim \log a_n$).

\begin{proof}[Proof of Proposition \ref{prop: general tightness}]

Let $P_n \subset \mathbb{R}^{n}$ be the polytopes constructed in Proposition \ref{prop: tightness with a} for $a = \frac{1}{2}$.
Let $N$ be given by $n =  \lfloor f(N) \rfloor^2$ (here we use $f(n)<\sqrt{n}$). Let $k_ N = k = n\cdot \lfloor \frac{N}{n} \rfloor$ (note that $k \sim N$). We orthogonally decompose $\mathbb{R}^k = \mathbb{R}^{n} \bigoplus .... \bigoplus \mathbb{R}^{n}$, where there are $\lfloor \frac{N}{n} \rfloor$ summands. Let $P_k = P_n \times ... \times P_n \subset \mathbb{R}^k$. We have, for some absolute constants $c_, c'>0$:
    \[ |V_k| = |V_n|^{\frac{k}{n}} \sim 2^{\frac{ck}{f(k)}},  \]
    \[ |\mathcal{F}_k| = \frac{k}{n}|\mathcal{F}_n| \sim  \frac{k}{f(k)^2} 2^{c'f(k)}.  \]
So,
\[ \log|V_k| \cdot \log|\mathcal{F}_k| \sim \frac{k}{f(k)}\left( f(k) + \log\left(\frac{k}{f(k)^2}\right) \right) \sim k\left( 1 + \frac{\log k}{f(k)}\right) \]
Moreover, by choosing $f(N) = \log(\alpha)$ (for $\alpha > 3$, say) we can also saturate the refined FLM inequality (corollary \ref{cor: refined FLM}). We have
\[ |V_k| \sim 2^{\frac{ck}{\log \alpha}}, \quad |\mathcal{F}_k| \sim \frac{\alpha k}{\log^2\alpha}, \]
and so
\[ \log|V_k| \cdot \log|\mathcal{F}_k| \sim k \cdot \frac{\log(\alpha k)}{\log\alpha} \]
\end{proof}

\subsection{Tightness of the general FLM inequality}

While the question of tightness of the FLM inequality has been essentially settled in the previous section, very little is known in every other case. To start, one can ask, in light of Proposition \ref{prop: tightness with a}, if given $a,b,c$ such that $a+b+2c=2$, one can find polytopes $P_n \subset \mathbb{R}^n$ with 
\[ \log|V_n| \sim n^a, \quad \log|\mathcal{F}_n| \sim n^b, \:\: and \: \: \left(\frac{R}{r}\right) \sim n^c \: ? \] 
\noindent Proposition \ref{prop: tightness with a} showed that for $c=\frac{1}{2}$, the answer is yes for all $a$ (so that $b=1-a$). Let us start by explaining some obvious restrictions that need to be imposed on $a,b,c$ for this to generally hold. The following theorem is classical, and can be found in \cite[Section 3]{AGA}.

\begin{thm}\label{thm: keith ball bound}
    Let $P \subset \mathbb{R}^n$ be a polytope such that $rB_2 \subset P \subset RB_2^n$. Then \\ $\log|V| \gtrsim n \left(\frac{r}{R} \right)^2 = n^{1-2c}$, and by duality, also $\log|\mathcal{F}| \gtrsim n^{1-2c}$.
\end{thm}

\noindent Now it follows that for every $c\in (0,1)$, we must have $a,b \geq 1-2c$. Note that this condition is non-trivial for $c < \frac{1}{2}$. We conjecture that this necessary condition is also sufficient.

\begin{conj}\label{conj: FLM}
    For all numbers $a,b,c >0$ such that $c < 1$, $a+b+2c=2$ and $a,b \geq 1-2c$, the following holds: There exists a sequence of polytopes $P_n \subset \mathbb{R}^n$ with
    \[ \log|V_n| \sim n^a, \:\: \log|\mathcal{F}_n| \sim n^b \:\: and \:\: \left( \frac{R_n}{r_n} \right)^2 \sim n^{2c} . \]
\end{conj}

\noindent We finish this section with some examples of constructions which saturate the FLM inequality in some special cases. 

\subsubsection{Examples of polytopes that saturate (or nearly saturate) the general FLM inequality}

The following theorem will be useful in the examples that follow. It is called the \textsl{low-$M^{*}$ estimate} and is due to Milman; see \cite{Mil}, and also \cite[Chapter 7]{AGA} for a detailed exposition and the optimal version we use here:

\begin{thm}\label{thm: low M*}(Low-$M^{\ast}$ estimate)
Let $K \subset \R^{n}$ be a convex body with the origin in its interior.  
Let $1 \le k \le n$ be an integer and write $k=\lambda n$.  
Then a random subspace $E \in G_{n,k}$ satisfies
\[
  \diam\bigl(K \cap E\bigr)
  \lesssim \frac{1}{\sqrt{1-\lambda}} M^{\ast}(K)
\]
with probability at least $1 - e^{-\tilde{c}\,(1-\lambda)n}$.
\end{thm}

\noindent We state the first example as a proposition:
\noindent \begin{prop}\label{prop: random section 1} Conjecture \ref{conj: FLM} Holds for all triplets $(a,b,c) = (1,b,\frac{1-b}{2})$, for all $b \in (0,1)$ (and by duality, it holds for all triplets $(a,b,c) = (b,1,\frac{1-b}{2})$ as well).
    
\end{prop}
\begin{proof} Given $a \in (0,1)$, let $P_n^a$ be the polytopes from \ref{prop: tightness with a}, satisfying $\log|V_n| \sim n^a$ and $\log|\mathcal{F}_n| \sim n^{1-a}$. We scale $P_n^a$ so that $B_2^n \subset P_n^a \subset \sqrt{n} B_2^n$. 

\noindent Let $E \subset \mathbb{R}^{2n}$ be a random subspace of dimension $n$. Let $Q_n^a = P_{2n}^a \cap E \subset \mathbb{R}^n$. We have:
\begin{itemize}
    \item $\log|\mathcal{F}(Q_n^a)| \lesssim n^{1-a}$, as the number of facets can only decrease by sections.
    \item $\log|V(Q_n^a)| \lesssim n$, as vertices of an $n$-dimensional section (generically) can only arise from $n$ dimensional faces of $P_{2n}^a$, and there are only exponentially many such faces, as this is a general (and easy to prove) property of Hanner polytopes.
    \item $R(Q_n^a) \lesssim M^*(P_{2n}^a) \lesssim n^{\frac{a}{2}}$ with high probability, by Theorem \ref{thm: low M*} and lemma \ref{lemma: mean width bound}.
\end{itemize}  

\noindent Thus, with non-zero probability, we get that the random section $E \cap P_{2n}^{a} \subset \mathbb{R}^n$ saturates FLM with the parameters stated. Note that all the implied constants above depend only on $a$.
\end{proof}

\noindent \textbf{Remark.} The following argument does not saturate FLM but can come close, and so may be of interest as well. As in proposition \ref{prop: random section 1}, we take a random section of $P_n^a$, but this time of dimension $k := \lfloor n - n^{\delta} \rfloor$ for some fixed $\delta \in (0,1)$. Let $Q_n^{a,\delta} = P_n^a \cap E$ where $E \subset \mathbb{R}^n$ is a random subspace of dimension $k$. We have:
\begin{itemize}
    \item $\log|\mathcal{F}(Q_n^{a,\delta})| \lesssim n^{1-a}$.
    \item $\log|V(Q_n^{a,\delta})| \lesssim \log\binom{2^{n^a}}{n^\delta} \lesssim n^{a+\delta}$, as vertices of an $n-n^{\delta}$-dimensional section (generically) can only arise from $n^{\delta}$ dimensional faces of $P_{2n}^a$.
    \item $R(Q_n^{a,\delta}) \lesssim n^{\frac{1-\delta}{2}}M^*(P_{2n}^a) \lesssim n^{\frac{1-\delta+a}{2}}$ with high probability, by Theorem \ref{thm: low M*} and lemma \ref{lemma: mean width bound}.
\end{itemize}  
So:
\[ \log|V| \cdot |\log|\mathcal{F}| \cdot \left( \frac{R}{r} \right)^2 \sim n^{2+a}. \]
Since the product is independent of $\delta$, we may apply this to all $\delta \in (0,1)$ and get `almost saturation' for many different parameters.
\\

\noindent \textbf{2.} The following example shows that the general FLM inequality is tight when $\frac{R}{r} < C$ for some universal constant $C$. In this case, the inequality becomes
\[ \log|V| \cdot \log|\mathcal{F}| \gtrsim n^2. \]
By Theorem \ref{thm: keith ball bound} we also have $\log|V|, \log|\mathcal{F}| \gtrsim n$, and so to saturate the inequality, we need polytopes $P_n$ for which $\frac{R_n}{r_n} < C$ and $\log|V_n| \sim \log|\mathcal{F}_n| \sim n$. This can be achieved as follows: Let $n \in \mathbb{N}$, and let $B_1^{2n} \subset \mathbb{R}^{2n}$ be the $2n$-dimensional cross-polytope. Let $E$ be a random $n$-dimensional subspace of $\mathbb{R}^{2n}$. The polytope $P_n := B_1^{2n} \cap E \subset \mathbb{R}^n$ satisfies $|\mathcal{F}_n| \leq 2^{2n}$ (as sections of a polytope can only decrease its number of facets) and $|V_n| \lesssim c^n$ for some $c>0$ (as a vertex of $P_n$ is given by the intersection of $E$ with an $n$ dimensional face of $B_1^{2n}$, and there are only exponentially many such faces). Finally, by Milman's version of Dvoretzky's theorem (see \cite[Section 5]{AGA} for an exposition) we have that $cB_2^n \subset E \cap B_1^{2n}$ for some constant $c>0$, with high probability (note that always $E \cap B_1^{2n} \subset B_2^n$). Thus, an $n$-dimensionsal random section of $B_1^{2n}$ saturates the FLM inequality in this case, with high probability on choosing the section (all that matters for us is that this probability is non-zero).
\\

\noindent \textbf{3.} Let $S_n \subset \mathbb{R}^n$ be a \textbf{regular} simplex satisfying $B_2^n \subset S_n \subset nB_2^n$. Let $k = n - f(n)$ for some rate function $f(n)$ for which $f(n) = o(n)$ (for example, $f(n) = n^{\delta}$ for some $\delta \in (0,1)$). Let $E$ be a random, uniformly chosen $k$-dimensional subspace of $\mathbb{R}^n$. We have (this holds for all subspaces of dimension $k$)
\[ |\mathcal{F}| \leq n+1  \:\: and \:\: \log |V| \lesssim \log \binom{n}{\lfloor f(n) \rfloor} \lesssim f(n)\left( \log n +1\right). \]
Moreover, it is standard to verify that $M^{*}(S_n) \sim \sqrt{n\log n}$ (see \cite[Table 4.1]{SA}), and so by Theorem \ref{thm: low M*} we have, with high probability:
\[ R(S \cap E) \lesssim \sqrt{\frac{n}{f(n)}}\sqrt{n\log n} = n \sqrt{\frac{\log n}{f(n)}}. \]
And so (note that trivially $r \geq 1$),
\[ \log|V| \cdot \log|\mathcal{F}| \cdot \left(\frac{R}{r} \right)^2 \sim \log n \cdot f(n)\log n \cdot n^2 \frac{\log n}{f(n)} = n^2\log^3 n . \]

\noindent Choosing $f(n) = n^{\delta}$, for example, gives $\log|V| \sim n^{\delta} \log n$ and $\left( \frac{R}{r} \right)^2 \sim n^{2-\delta} \log n$. Note that $|\mathcal{F}| = n(1+o(1))$. This saturates FLM up to a $\log^3 n$ factor (actually up to a $\log^2 n$ factor, by (2)).
\\

\noindent \textbf{4.} The following is not an example, but rather an attempt to highlight a special case of the FLM inequality. Setting $\frac{R}{r} \sim \frac{n}{\log n}$ gives the inequality
\[ \log|V| \cdot \log|\mathcal{F}| \gtrsim \log^2n, \]
which is, of course, trivial. However, I do not know whether it can be saturated under the above assumption on $\frac{R}{r}$. Note that for the regular simplex $S_n$ from the previous example, we have $\frac{R}{r} = n$ and $|V| = |\mathcal{F}| = n+1$; so in other words, we are asking if we can make $\frac{R}{r}$ smaller by a $\log$ factor, while paying the `small price' of making $|V|$ and $ |\mathcal{F}|$ grow polynomially rather than linearly. As far as I am aware, this is not known.
\\

\section{A geometric version of FLM}
This final section is meant to highlight an inequality which is formally weaker than the FLM inequality, and does not involve polytopes, but rather the \textbf{Dvoretzky dimension} of convex bodies. The following was proved by Milman (see \cite[Chapter 5]{AGA} for an exposition):
\begin{thm}
    Let $K \subset \mathbb{R}^n$ be a convex body with $rB_2^n \subset K \subset RB_2^n$. There exists a subspace $E$ of dimension $k \sim \eps^2n(M(K)r)^2$ and an ellipsoid $\mathcal{E} \subset E$ such that:
    \[ (1-\eps)\mathcal{E} \subset K\cap E \subset (1+\eps)\mathcal{E}. \]
    By duality, there is an orthogonal projection to a subspace $E$ of dimension $k \sim \eps^2 n \left(\frac{M^{*}(K)}{R}\right)^2$ and an ellipsoid $\mathcal{E} \subset E$ such that
    \[ (1-\eps)\mathcal{E} \subset Proj_E(K) \subset (1+\eps)\mathcal{E}.  \]

\end{thm}

\noindent Let $dv_S(K) := n(M(K)r)^2$ and $dv_P(K) = n\left(\frac{M^{*}(K)}{R}\right)^2$ denote the Dvoretzky intersection dimension and Dvoretzky projection dimension, respectively. Clearly, we have:
\begin{equation} dv_S \cdot dv_P \cdot \left(\frac{R}{r}\right)^2 = n^2(MM^{*})^2 \geq n^2  .
\end{equation}
Moreover, by a similar argument to the one presented in the proof of corollary $\ref{cor: refined FLM},$
\[ \log|V| \gtrsim dv_P, \quad \log|\mathcal{F}| \gtrsim dv_S .\]
Thus we can conclude that (4) is formally a weaker inequality than the FLM inequality. Saturating (4) is no longer a question about polytopes; one needs to find convex bodies with prescribed Dvoretzky dimensions. And indeed, we can, in the following sense:
\begin{prop}\label{prop: geometric FLM}
    For all numbers $a,b,c < 1$ such that $a+b+2c=2$ and $a,b \geq 1-2c$, the following holds: There exists a sequence of convex bodies $K_n \subset \mathbb{R}^n$ with
    \[ dv_P(K_n) \sim n^a, \:\: dv_S(K_n)\sim n^b \:\: and \:\: \left( \frac{R_n}{r_n} \right)^2 \sim n^{2c} . \]
\end{prop}
\noindent \textbf{Remarks.}    1. Note the restriction $a,b < 1$, which doesn't appear in conjecture \ref{conj: FLM}. This comes from the trivial inequalities $dv_S, dv_P \leq n$.
Also, note that the restriction $a,b \geq 1-2c$, originated in Theorem \ref{thm: keith ball bound}, is needed in this variant as well; indeed, it is standard to verify that if $\frac{R}{r} \lesssim n^c$, then $dv_S, dv_P \gtrsim n^{1-2c}$.

\noindent 2. Corollary \ref{cor: refined FLM} demonstrates that the `refined FLM' inequality in fact holds in the form of (3). Indeed, if $|\mathcal{F}| = \alpha n$ and $\frac{R}{r} \lesssim \sqrt{n}$, we have $dv_{P} \gtrsim \sqrt{\frac{n}{\log(1+\alpha)}}$, and we always have $dv_{S} \gtrsim \log n$ (again by \cite[Chapter 5]{AGA}). So, 
\[ dv_P \cdot dv_S \gtrsim n \cdot \frac{\log (n)}{\log(1+\alpha)} .\]
\newpage
We now turn to prove Proposition \ref{prop: geometric FLM}. For all numbers $\beta, c$ such that $c \in (0,1)$ and $\beta \in (c-\frac{1}{2}, c)$, let $K_n^{c,\beta} \subset \mathbb{R}^n$ be the following convex body:
\[ K_n^{c,\beta} := \conv \left( n^{\beta} B_2^n \cap \{|\langle x, e_1 \rangle | \leq 1 \}, n^c e_n \right)
 \]
 It is clear that we have $B_2^n \subset K_n^{c, \beta} \subset n^c B_2^n$ and that both of the inclusions are tight.
 \begin{lem}
     for all $c,\beta$ as above we have
     \[ M^{*}(K_n^{c,\beta}) \sim n^{\beta} \quad and \quad M(K_n^{c,\beta)}) \sim n^{-\beta}. \]
 \end{lem}
\begin{proof}
    To bound from above, we use the fact that $\conv(A, B) \subset A+B$ for any convex sets $A,B$, and the fact that for a unit interval $I \subset \mathbb{R}^n$ we have $M^{*}(I) \sim \frac{1}{\sqrt{n}}$:
    \[ M^{*}(K_n^{c, \beta}) \leq M^{*}(n^\beta B_2^n + n^c e_n) \sim n^\beta + n^{c-\frac{1}{2}} \sim n^{\beta}.  \]
    For the other bound, let $A_\beta = \{ x \in S^{n-1}: |\langle x, e_1 \rangle| \leq n^{-\beta}\}$. By lemma \ref{lemma: concentration of measure}:
    \[ \sigma(A_\beta) \geq 1-e^{-\frac{1}{2}n^{1-2\beta}} \]
    The essential observation regarding $A_\beta$ is that for all $x \in A_\beta$ we have \[ h_{K_n^{c,\beta}}(x) = h_{\conv(n^{\beta}B_2^n, n^c e_n)}(x), \]
    and so:
    \[ M^{*}(K_n^{c, \beta}) \geq \int_{A_\beta}h_{K_n^{c, \beta}}d\sigma \geq n^{\beta}(1-e^{-\frac{1}{2}n^{1-2\beta}}) \sim n^{\beta}, \]
    where in the last step we used $\beta < \frac{1}{2}$. 
    For the estimate of $M(K_{n}^{c,\beta})$: from the inequality $MM^{*} \geq 1$ we have $ M(K_n^{c, \beta}) \gtrsim n^{-\beta}$, and for the other side:
    \[  M(K_n^{c, \beta}) = \int_{S^{n-1}}\|x\|_{K_n^{C,\beta}}d\sigma \leq \int_{S^{n-1}}\|x\|_{n^{\beta}B_2^n \cap \{|\langle x, e_1 \rangle | \leq 1\}}d\sigma  =   \int_{A_{\beta}}\|x\|_{n^{\beta}B_2^n \cap \{|\langle x, e_1 \rangle | \leq 1\}}d\sigma \: +  \]
    \[\int_{S^{n-1} \setminus A_\beta}\|x\|_{n^{\beta}B_2^n \cap \{|\langle x, e_1 \rangle | \leq 1\}}d\sigma \leq n^{-\beta} + \int_{S^{n-1} \setminus A_{\beta}}\|x\|_{n^{\beta}B_2^n \cap \{|\langle x, e_1 \rangle | \leq 1\}}d\sigma \leq n^{-\beta} + e^{-\frac{1}{2}n^{1-2\beta}} \sim n^{-\beta}. \]
    \end{proof}
    
\begin{proof}[Proof of Proposition \ref{prop: geometric FLM}]
Let $a,b,c$ as in the proposition be given. Let $\beta := c + \frac{a-1}{2}$. We have, for $K_n^{c,\beta}$:
\begin{itemize}
    \item $dv_P(K_n^{c,\beta}) \sim n^{1+2\beta-2c} = n^a $
    \item $dv_S(K_n^{c,\beta}) \sim n^{1-2\beta} = n^b $
    \item $\left( \frac{R_n^{c,\beta}}{r_n^{c,\beta}} \right)^2 = n^{2c}$
\end{itemize} 

\end{proof}

\noindent \textbf{Acknowledgments.} I would like to thank Guillaume Aubrun and Alexander Barvinok for useful discussions.

\bibliographystyle{amsplain}

\begin{thebibliography}{99}

\bibitem{AGA}
Artstein-Avidan S., Giannopoulos A., Milman V.D.,
{\em Asymptotic Geometric Analysis, Part I}. Mathematical Surveys and Monographs,  volume 202, American Mathematical Society, Providence, RI. (2015).
https://doi.org/10.1090/surv/202


\bibitem{FLM}
Figiel, T., Lindenstrauss, J.,  Milman, V.D. {\em The dimension of almost spherical sections of convex bodies}. Acta Math. 139, 53–94 (1977). https://doi.org/10.1007/BF02392234.

\bibitem{Bar}
Barvinok, A. {\em  A bound for the number of vertices of a polytope with applications}. Combinatorica 33, 1–10 (2013). https://doi.org/10.1007/s00493-013-2870-9.

\bibitem{SA}
Aubrun, G.,  Szarek, S. J. (2017). {\em Alice and Bob Meet Banach: The Interface of Asymptotic Geometric Analysis and Quantum Information Theory} (Vol. 223). American Mathematical Society.

\bibitem{Mil}
V. D. Milman, {\em Random subspaces of proportional dimension of finite-dimensional normed
spaces: approach through the isoperimetric inequality}, Banach spaces (Columbia, Mo.,
1984), Lecture Notes in Mathematics, no. 1166, Springer-Verlag, Berlin, 1985, pp. 106–
115. MR 87j:46037b

\bibitem{Glus}
Gluskin, E.D. {\em Extremal properties of orthogonal parallelepipeds and their applications to the geometry of Banach spaces}. Mathematics of the USSR‑Sbornik, 64(1), 85–96. https://doi.org/10.1070/SM1989v064n01ABEH003295


\end{thebibliography}

\end{document}